\documentclass{article}
\usepackage{graphicx} 
\usepackage[utf8]{inputenc}
\usepackage[margin=1in]{geometry}
\usepackage{amsmath}
\usepackage{amsthm}
\usepackage{amsfonts}
\usepackage{amssymb}
\usepackage{txfonts}
\usepackage{cite}
\usepackage{ragged2e}
\usepackage{authblk}
\usepackage{hyperref}

\providecommand{\keywords}[1]
{
  \small	
  \textbf{\textit{Keywords: }} #1
}

\providecommand{\subjclass}[1]
{
  \small	
  \textbf{\textit{Mathematics Subject Classification: }} #1
}

\newtheorem{theorem}{Proposition}
\newtheorem{corollary}[theorem]{Corollary}
\title{On a double series}
\author{\emph{Aung Phone Maw}}
\date{October 2025}

\begin{document}

\maketitle

\begin{abstract}
    We shall investigate and arrive at a certain functional property of the double series \[
 \sum\limits_{n,r\geq 1}\frac{1}{\sqrt{x^2n^2+r^2+w^2}\left( e^{2 \pi y\sqrt{x^2n^2+r^2+w^2}}-1\right)}.
\]
    
\end{abstract}
. \\
.

\section{Introduction}

Let us recall that by $(3.11)$ in \cite{apm} we have 
\begin{multline}
    \frac{1}{x_2} \sum\limits_{n=1}^{\infty}\frac{1}{\sqrt{x_1^2n^2+w^2}\left( e^{\frac{2 \pi}{x_2}\sqrt{x_1^2n^2+w^2}}-1\right)} - \frac{1}{x_1} \sum\limits_{n=1}^{\infty}\frac{1}{\sqrt{x_2^2n^2+w^2}\left( e^{\frac{2 \pi}{x_1}\sqrt{x_2^2n^2+w^2}}-1\right)} \\= \frac{1}{2wx_1\left(e^{\frac{2\pi w}{x_1}}-1\right)} - \frac{1}{2wx_2\left(e^{\frac{2\pi w}{x_2}}-1\right)} + \frac{1}{2} \sum\limits_{n=1}^{\infty} \left\{ \frac{1}{x_1\sqrt{x_2^2n^2+w^2}}- \frac{1}{x_2\sqrt{x_1^2n^2+w^2}}\right\} + \frac{\ln x_1}{2 x_1 x_2} -\frac{\ln x_2}{2 x_1 x_2} +\frac{1}{4wx_1}- \frac{1}{4wx_2}. \label{transform} \tag{1.1}
\end{multline}

Now let us define 
\[
f(x,y;w) = \sum\limits_{n,r \geq 1}\frac{1}{\sqrt{x^2n^2+r^2+w^2}\left( e^{2 \pi y\sqrt{x^2n^2+r^2+w^2}}-1\right)}, f(x,y;0)=f(x,y).
\]
We have the obvious relation 
\[
 x^{-1}f(x^{-1},y;wx^{-1}) = f(x,yx^{-1};w). \label{ftrans} \tag{1.2}
\]
Then we shall provide the following functional relation for $f$.

\begin{theorem} For complex values $x_1, x_2,$ and $w$, there holds 
\begin{flalign*}
  x_1f(x_1,x_2^{-1};w)-x_2f(x_2,x_1^{-1};w) -f(x_1,x_1^{-1}x_2^{-1};w)+x_2f(x_2,x_1;wx_1^{-1}) \\ -x_1f(x_1,x_2;wx_2^{-1})+f(x_2,x_1^{-1}x_2^{-1};w) +f(x_1,x_2x_1^{-1};wx_2^{-1})-f(x_2,x_1x_2^{-1};wx_1^{-1}) \\ = \frac{1}{4}\left(x_1 - \frac{1}{x_1} \right) \ln x_2 -\frac{1}{4}\left( x_2 - \frac{1}{x_2}\right) \ln x_1 \\+ \frac{1}{2} \sum\limits_{r=1}^{\infty}\sum\limits_{n=1}^{\infty}  \left\{ \frac{x_2}{\sqrt{x_2^2n^2+r^2+w^2}}- \frac{x_1}{\sqrt{x_1^2n^2+r^2+w^2}}  -\frac{x_1x_2}{\sqrt{x_2^2x_1^2 n^2+x_1^2r^2+w^2}}+\frac{1}{\sqrt{n^2+r^2x_1^2+w^2}} \right. \\ \left. -\frac{1}{\sqrt{n^2+x_2^2r^2+w^2}}+ \frac{x_1x_2}{\sqrt{x_1^2x_2^2n^2+x_2^2r^2+w^2}}  +\frac{x_1}{\sqrt{x_1^2 n^2+x_1^2x_2^2r^2+w^2}}-\frac{x_2}{\sqrt{x_2^2n^2+r^2x_1^2x_2^2+w^2}} \right\}, \label{proposition1} \tag{1.3}
\end{flalign*}
provided that all the series are convergent.
\end{theorem}
For natural numbers $a,b$ and $N$ let $r_{a,b}(N)$ be defined as the number of ways $N$ can be written in the form $ax^2+by^2$ where $x$ and $y$ are natural numbers. Then we shall also provide the following corollary of \textbf{Proposition 1}. 

\begin{corollary}
    For natural numbers $a,b,c,$ and $d$, there holds
    \begin{flalign*}
        \sum\limits_{n=1}^{\infty} \frac{r_{a,b}(n)}{\sqrt{n}}\left\{ \frac{\sqrt{a}}{e^{2\pi\sqrt{ncb^{-1}d^{-1}}}-1}-\frac{\sqrt{a}}{e^{2\pi\sqrt{ndc^{-1}b^{-1}}}-1}+\frac{\sqrt{b}}{e^{2\pi\sqrt{ndbc^{-1}a^{-1}}}-1}- \frac{\sqrt{b}}{e^{2\pi\sqrt{ncbd^{-1}a^{-1}}}-1}\right\} \\ -  \sum\limits_{n=1}^{\infty} \frac{r_{c,d}(n)}{\sqrt{n}}\left\{ \frac{\sqrt{c}}{e^{2\pi\sqrt{nab^{-1}d^{-1}}}-1}-\frac{\sqrt{c}}{e^{2\pi\sqrt{nba^{-1}d^{-1}}}-1}+\frac{\sqrt{d}}{e^{2\pi\sqrt{ndbc^{-1}a^{-1}}}-1}- \frac{\sqrt{d}}{e^{2\pi\sqrt{nadb^{-1}c^{-1}}}-1}\right\} \\= \frac{\left( \sqrt{\frac{c}{d}}- \sqrt{\frac{d}{c}}\right) \ln\left( \frac{a}{b}\right)}{8}-\frac{\left( \sqrt{\frac{a}{b}}- \sqrt{\frac{b}{a}}\right) \ln\left( \frac{c}{d}\right)}{8}. \label{sumofsquaresseries} \tag{1.4}
    \end{flalign*}
\end{corollary}

\section{Proof of Propsition 1}

In \eqref{transform} substitute $\sqrt{r^2+w^2}$ for $w$ and sum over all natural numbers $r$ to get

\begin{flalign*}
    \frac{1}{x_2}f(x_1,x_2^{-1};w)-\frac{1}{x_1}f(x_2,x_1^{-1};w) = \frac{1}{2x_1}\sum\limits_{r} \frac{1}{\sqrt{r^2+w^2}\left( e^{\frac{2\pi}{x_1}\sqrt{r^2+w^2}}-1\right)}- \frac{1}{2x_2}\sum\limits_{r} \frac{1}{\sqrt{r^2+w^2}\left( e^{\frac{2\pi}{x_2}\sqrt{r^2+w^2}}-1\right)} \\ + \sum\limits_{r=1}^{\infty} \left\{ \frac{1}{2} \sum\limits_{n=1}^{\infty} \left\{ \frac{1}{x_1\sqrt{x_2^2n^2+r^2+w^2}}- \frac{1}{x_2\sqrt{x_1^2n^2+r^2+w^2}}\right\} + \frac{\ln x_1}{2 x_1 x_2} -\frac{\ln x_2}{2 x_1 x_2} +\frac{1}{4x_1\sqrt{r^2+w^2}}- \frac{1}{4x_2\sqrt{r^2+w^2}} \right\}. \label{eq1} \tag{2.1}
\end{flalign*}
Multiply with $x_1$ to both sides, then substitute $x_1^{-1}$ for $x_1$, $wx_1^{-1}$ for $w$, and divide both sides by $x_1$ to arrive at

\begin{flalign*}
    \frac{1}{x_2x_1^2}f(x_1^{-1},x_2^{-1};wx_1^{-1})-\frac{1}{x_1}f(x_2,x_1;wx_1^{-1}) = \frac{1}{2}\sum\limits_{r} \frac{1}{\sqrt{r^2x_1^2+w^2}\left( e^{2\pi\sqrt{r^2x_1^2+w^2}}-1\right)}- \frac{1}{2x_2x_1}\sum\limits_{r} \frac{1}{\sqrt{r^2x_1^2+w^2}\left( e^{\frac{2\pi}{x_2x_1}\sqrt{r^2x_1+w^2}}-1\right)} \\ + \sum\limits_{r=1}^{\infty} \left\{ \frac{1}{2} \sum\limits_{n=1}^{\infty} \left\{ \frac{1}{\sqrt{x_2^2x_1^2 n^2+x_1^2r^2+w^2}}- \frac{1}{x_2x_1\sqrt{n^2+r^2x_1^2+w^2}}\right\} - \frac{\ln x_1}{2 x_1 x_2} -\frac{\ln x_2}{2 x_1 x_2} +\frac{1}{4\sqrt{r^2x_1^2+w^2}}- \frac{1}{4x_2x_1\sqrt{r^2x_1^2+w^2}} \right\}. \label{eq2} \tag{2.2}
\end{flalign*}
Now subtract \eqref{eq2} from \eqref{eq1} and use \eqref{transform} to get

\begin{flalign*}
    \frac{1}{x_2}f(x_1,x_2^{-1};w)-\frac{1}{x_1}f(x_2,x_1^{-1};w) -\frac{1}{x_2x_1^2}f(x_1^{-1},x_2^{-1};wx_1^{-1})+\frac{1}{x_1}f(x_2,x_1;wx_1^{-1}) = \frac{1}{4w\left( e^{2\pi w}-1\right)}-\frac{1}{4wx_1\left(e^{\frac{2\pi w}{x_1}}-1\right) } \\-\frac{\ln x_1}{4x_1} + \frac{1}{8w}-\frac{1}{8wx_1}- \frac{1}{2x_2}\sum\limits_{r} \frac{1}{\sqrt{r^2+w^2}\left( e^{\frac{2\pi}{x_2}\sqrt{r^2+w^2}}-1\right)} +\frac{1}{2x_2x_1}\sum\limits_{r} \frac{1}{\sqrt{r^2x_1^2+w^2}\left( e^{\frac{2\pi}{x_2x_1}\sqrt{r^2x_1+w^2}}-1\right)} \\ + \sum\limits_{r=1}^{\infty} \left\{ \frac{1}{2} \sum\limits_{n=1}^{\infty} \left\{ \frac{1}{x_1\sqrt{x_2^2n^2+r^2+w^2}}- \frac{1}{x_2\sqrt{x_1^2n^2+r^2+w^2}}  -\frac{1}{\sqrt{x_2^2x_1^2 n^2+x_1^2r^2+w^2}}+\frac{1}{x_2x_1\sqrt{n^2+r^2x_1^2+w^2}}\right\} \right. \\ \left. + \frac{\ln x_1}{ x_1 x_2}  - \frac{1}{4x_2\sqrt{r^2+w^2}} + \frac{1}{4x_2x_1\sqrt{r^2x_1^2+w^2}}\right\}. \label{eq3} \tag{2.3}
\end{flalign*}
Multiply both sides with $x_2$, substitute $x_2^{-1}$ for $x_2$,  $ wx_2^{-1}$ for $w$, then divide by $x_2$ to yield

\begin{flalign*}
    \frac{1}{x_2}f(x_1,x_2;wx_2^{-1})-\frac{1}{x_1x_2^2}f(x_2^{-1},x_1^{-1};wx_2^{-1}) -\frac{1}{x_2x_1^2}f(x_1^{-1},x_2;wx_2^{-1}x_1^{-1})+\frac{1}{x_1x_2^2}f(x_2^{-1},x_1;wx_2^{-1}x_1^{-1}) \\
    = \frac{1}{4x_2w\left( e^\frac{2\pi w}{x_2}-1\right)}  -\frac{1}{4wx_2x_1\left(e^{\frac{2\pi w}{x_2x_1}}-1\right) } -\frac{\ln x_1}{4x_1x_2^2} + \frac{1}{8wx_2}-\frac{1}{8wx_2x_1}  \\
        - \frac{1}{2}\sum\limits_{r} \frac{1}{\sqrt{x_2^2r^2+w^2}\left( e^{2\pi\sqrt{x_2^2r^2+w^2}}-1\right)} +\frac{1}{2x_1}\sum\limits_{r} \frac{1}{\sqrt{r^2x_1^2 x_2^2+w^2}\left( e^{\frac{2\pi}{x_1}\sqrt{r^2x_1x_2^2+w^2}}-1\right)} \\ + \sum\limits_{r=1}^{\infty} \left\{ \frac{1}{2} \sum\limits_{n=1}^{\infty} \left\{ \frac{1}{x_1x_2\sqrt{n^2+x_2^2r^2+w^2}}- \frac{1}{\sqrt{x_1^2x_2^2n^2+x_2^2r^2+w^2}}  -\frac{1}{x_2\sqrt{x_1^2 n^2+x_1^2x_2^2r^2+w^2}}+\frac{1}{x_1\sqrt{x_2^2n^2+r^2x_1^2x_2^2+w^2}}\right\} \right. \\ \left. + \frac{\ln x_1}{ x_1 x_2}  - \frac{1}{4\sqrt{r^2x_2^2+w^2}} + \frac{1}{4x_1\sqrt{r^2x_1^2x_2^2+w^2}}\right\}. \label{eq4} \tag{2.4}
\end{flalign*}
Subtract \eqref{eq4} from \eqref{eq3} and use \eqref{transform} to get

\begin{flalign*}
    \frac{1}{x_2}f(x_1,x_2^{-1};w)-\frac{1}{x_1}f(x_2,x_1^{-1};w) -\frac{1}{x_2x_1^2}f(x_1^{-1},x_2^{-1};wx_1^{-1})+\frac{1}{x_1}f(x_2,x_1;wx_1^{-1}) \\ -\frac{1}{x_2}f(x_1,x_2;wx_2^{-1})+\frac{1}{x_1x_2^2}f(x_2^{-1},x_1^{-1};wx_2^{-1}) +\frac{1}{x_2x_1^2}f(x_1^{-1},x_2;wx_2^{-1}x_1^{-1})-\frac{1}{x_1x_2^2}f(x_2^{-1},x_1;wx_2^{-1}x_1^{-1}) \\ = \frac{1}{4w\left( e^{2\pi w}-1\right)}-\frac{1}{4wx_1\left(e^{\frac{2\pi w}{x_1}}-1\right) }-\frac{1}{4x_2w\left( e^\frac{2\pi w}{x_2}-1\right)}  +\frac{1}{4wx_2x_1\left(e^{\frac{2\pi w}{x_2x_1}}-1\right) }\\+ \frac{1}{4wx_2\left( e^{\frac{2\pi w}{x_2}}-1\right)}-\frac{1}{4w\left(e^{2\pi w}-1\right) }+ \frac{1}{4wx_1\left( e^{\frac{2\pi w}{x_1}}-1\right)}-\frac{1}{4wx_2x_1\left( e^{\frac{2\pi w}{x_2x_1}}-1\right)} \\ - \frac{\ln x_1}{4x_1}+\frac{\ln x_1}{4x_1x_2^2}+ \frac{\ln x_2}{4x_2}+ \frac{\ln x_1}{4x_1^2x_2} - \frac{\ln x_1x_2}{4x_1^2x_2}+ \frac{1}{8w}- \frac{1}{8wx_1}-\frac{1}{8wx_2}+ \frac{1}{8wx_2x_1}+\frac{1}{8wx_2}- \frac{1}{8w}+ \frac{1}{8wx_1}- \frac{1}{8wx_1x_2} \\ + \frac{1}{4} \sum\limits_{n=1}^{\infty} \left\{ \frac{1}{x_2\sqrt{n^2+w^2}}- \frac{1}{\sqrt{x_2^2n^2+w^2}}\right\} + \frac{1}{4} \sum\limits_{n=1}^{\infty} \left\{ \frac{1}{x_1\sqrt{x_1x_2^2n^2+w^2}}- \frac{1}{x_1x_2\sqrt{x_1^2n^2+w^2}}\right\} \\ + \sum\limits_{r=1}^{\infty} \left\{ \frac{1}{2} \sum\limits_{n=1}^{\infty} \left\{ \frac{1}{x_1\sqrt{x_2^2n^2+r^2+w^2}}- \frac{1}{x_2\sqrt{x_1^2n^2+r^2+w^2}}  -\frac{1}{\sqrt{x_2^2x_1^2 n^2+x_1^2r^2+w^2}}+\frac{1}{x_2x_1\sqrt{n^2+r^2x_1^2+w^2}}\right. \right. \\ \left.\left.-\frac{1}{x_1x_2\sqrt{n^2+x_2^2r^2+w^2}}+ \frac{1}{\sqrt{x_1^2x_2^2n^2+x_2^2r^2+w^2}}  +\frac{1}{x_2\sqrt{x_1^2 n^2+x_1^2x_2^2r^2+w^2}}-\frac{1}{x_1\sqrt{x_2^2n^2+r^2x_1^2x_2^2+w^2}} \right\} \right. \\ \left.  - \frac{1}{4x_2\sqrt{r^2+w^2}} + \frac{1}{4x_2x_1\sqrt{r^2x_1^2+w^2}} + \frac{1}{4\sqrt{r^2x_2^2+w^2}} -\frac{1}{4x_1\sqrt{r^2x_1^2x_2^2+w^2}}\right\}.
\end{flalign*}
 Multiply both sides by $x_1x_2$, and after some manipulations we finally arrive at \textbf{Proposition 1}. Now
put $w=0$ in \textbf{Proposition 1} to get 

\begin{corollary} For complex values $x_1,x_2$, there holds
\begin{flalign*} 
    x_2f(x_2,x_1)-x_1f(x_1,x_2)+x_1f(x_1, x_2^{-1}) - x_2 f(x_2, x_1^{-1}) + f(x_2, x_1^{-1}x_2^{-1}) - f(x_1, x_1^{-1}x_2^{-1})+ f(x_1, x_2 x_1^{-1}) - f(x_2,x_1x_2^{-1}) \\= \frac{1}{4}\left(x_1-\frac{1}{x_1}\right)\ln x_2 - \frac{1}{4}\left(x_2 - \frac{1}{x_2}\right) \ln x_1 , \label{cor3} \tag{2.5}
\end{flalign*}
provided that all the series are convergent.
\end{corollary}
Now put $x_2=\sqrt{\frac{a}{b}}$ and $x_1=\sqrt{\frac{c}{d}}$, where $a,b,c,d$ are natural numbers in \textbf{Corollary 3}, then we arrive at \eqref{sumofsquaresseries}.

\end{document}